\documentclass{amsart}

\usepackage{hyperref, color}
\usepackage{amssymb}

\theoremstyle{plain}
\newtheorem{theorem}                {Theorem}      [section]
\newtheorem*{theorem*}                {Theorem \ref{thm:appl}}

\theoremstyle{definition}

\numberwithin{equation}{section}

\newcommand{\field}[1]{\mathbb{#1}}
\newcommand{\real}{\field{R}}

\newcommand{\vep}{\varepsilon}

\begin{document}

\title
[The halfspace theorem for minimal hypersurfaces]
{The halfspace theorem for minimal hypersurfaces in regions bounded by minimal cones}

\author[Cavalcante]{Marcos Petr\'ucio Cavalcante}     
\address{Institute of Mathematics, Federal University of Alagoas, CEP 57072-970,
Macei\'o,  Alagoas,  Brazil}
\email{marcos@pos.mat.ufal.br}

\author[Costa-Filho]{Wagner Oliveira Costa-Filho}
\address{Campus Arapiraca, Federal University of Alagoas, CEP 57309-005, 
Arapiraca,  Alagoas, Brazil}
\email{fcow@bol.com.br}

\dedicatory{In memory of Professor Manfredo Perdig\~ao do Carmo}

\date{\today}

\begin{abstract}  
We prove that there are no minimal hypersurfaces properly immersed in any 
region of the Euclidean space bounded by unstable minimal  cones. 
{We also prove the analogous result for $r$-minimal hypersurfaces.}
\end{abstract}

\maketitle

\section{Introduction}

A classical problem on differential geometry is to determine the existence or non-existence 
of hypersurfaces immersed in certain regions of a given Riemannian manifold. 
In this context, important results have been obtained in the class of minimal immersions,
mainly motivated by the Calabi-Yau conjecture \cite{Calabi, Chern, Yau}, which asserts 
that a complete minimal hypersurface in $\mathbb R^{n+1}$ is unbounded. 
Many important results were done about this conjecture, including counterexamples
in its general form (see 
\cite{Alarcon}, 
%\cite{BJM}, 
\cite{ColdingMinicozzi}, 
\cite{JorgeXavier}, 
\cite{LopesMartinMorales}, 
\cite{MartinMorales},
\cite{Nadirashvili}, 
and  references therein). 

Related to this problem we have the celebrated halfspace theorem of Hoffman and Meeks 
\cite{HoffmanMeeks} in the Euclidean three-space. It says that  
\emph{there is no  proper, possibly branched, nonplanar minimal surface
contained in a halfspace of $\mathbb R^3$}.
This is a nice and clever application of the Maximum Principle for
minimal surfaces making use of the behavior of the catenoid in $\mathbb R^3$.
It is astonishing that this theorem does not hold in
higher dimensions and the catenoids $\mathcal C^n$ itself are simple counterexamples,
since they are contained in a slab given by two parallel hyperplanes 
of $\real^{n+1}$, when $n\geq 3$. 
%In fact, Rosenberg, Schulze and Spruck noticed 

It is then natural to investigate some sort of the halfspace theorem in higher
dimensions or in others contexts. In fact, there exits an active research on this
topic  and many results have been obtained. See for instance
\cite{Albanese,
Atsuji,
Bessa,
BJO, 
Bombieri, 
CE, 
Daniel, 
Impera,
Lopez, 
MariRigoli, 
MazetWanderley, 
Mazet1, 
Mazet2, 
Omori,
Rodriguez, 
RosenbergSchulzeSpruck}.

We point out that the first result about the non existence of minimal immersions 
in cones was done
by Omori in \cite{Omori}, where he considered cones strictly contained in the half-space. 
A very nice extension was latter obtained by Atsuji in \cite{Atsuji} using stochastic geometry.
We also draw attention to two interesting theorems regarding proper minimal hypersurfaces 
immersed in the Euclidean space $\mathbb R^{n+1}$. In  \cite{Lopez}, L\'opez proved
a nonimmersibility result considering cones over convex domains and in 
\cite{MariRigoli} Mari and Rigoli considered this problem for non-degenerated cones 
(see  also \cite{Ran}). These results were extended to $r$-minimal hypersurfaces in 
\cite{Albanese}.

In this note, we prove that the half-space property is valid when we consider 
unstable minimal cones as barriers. 
Before we state our theorem in a precise form we need to recall some definitions. 
Let $\Sigma\subset\mathbb S^{n}$ be a compact orientable  hypersurface embedded in the 
unit round sphere. 
Let us denote by $C_\Sigma\subset \real^{n+1}$ the closed half-cone 
over $\Sigma$,  that is, the mapping
of $\Sigma\times [0,\infty)$ given by $(m,t)\to tm$. It is simple to see that 
$C_\Sigma^*=C_\Sigma\setminus\{0\}$ is a minimal
hypersurface if and only if $\Sigma\subset \mathbb S^n$ is minimal. 
Moreover,  since $\Sigma$ is embedded, we have that 
$\mathbb R^{n+1}\setminus C_\Sigma$ has two connected components. 
In \cite{Simons}, Simons proved that if 
$\Sigma\subset \mathbb S^n$ is not totally geodesic and $3\leq n \leq 6$, 
then $C_\Sigma\subset \real^{n+1}$ is an {unstable} minimal hypersurface 
(see also \cite{SSY}). 
%Moreover, Simons gave an example showing that it fails in higher dimensions, $n\geq 7.$
We recall that a minimal hypersurface is said to be \emph{stable} if it minimizes the area up
second order and \emph{unstable} otherwise. 
Since  the second variation of the area functional of a minimal hypersurface $M$ is
given by the Jacobi operator $J=\Delta + \|A\|^2$, where $A$ is the second fundamental
form of $M$, we conclude that $M$ is unstable if and only if 
$J$ has some negative eigenvalue on some compact domain. 
Our main theorem is the following.

\begin{theorem}\label{main}
Let $C_\Sigma \subset \mathbb R ^{n+1}$ be an unstable minimal cone, $n\geq 3$. Then
there is no minimal hypersurface $M^n$ properly immersed in $\mathbb R^{n+1}\setminus C_\Sigma$.
\end{theorem}

The proof is inspired by the proof of the half-space theorem for self-shrinkers by the first author and
Espinar \cite{CE} and the idea is to use the unstability of $C_\Sigma$ and the cone structure 
to create a perturbation of $C_\Sigma$ that lead a contradiction to the Maximum Principle.

In the following we denote by $\bar M(c)^{n+1}$ the space form of constant sectional
curvature $c$, which for simplicity, we assume to be $-1,0$ or $1$.
We recall that, given $r\in\{1,\ldots, n\}$, 
a hypersurface $M^n\subset\bar M(c)^{n+1}$ is said to be
\emph{$r$-minimal} if the symmetric function $S_{r}$ 
of the principal curvatures of $M$ is zero.
This is a natural generalization of minimal hypersurfaces ($r=1$) and also of
zero scalar curvature hypersurfaces ($r=2$) in the Euclidean space. 
In \cite{Reilly}, Reilly proved that such hypersurfaces are critical points of a certain  
$r$-area functional (see Section \ref{rminimal} for details) and computed its second variation 
formula (see also \cite{BC}). The second variation is given by a quadratic form in terms of a
second order linear operator, whose principal part is given in a divergent form (\cite{Rosenberg}). 
In this case, we say that a $r$-minimal hypersurface is \emph{$r$-stable} if its second variation
is nonnegative and \emph{$r$-unstable} otherwise (see \cite{BC, ACE}). 
As in the minimal case, a direct computation shows that $\Sigma^{n-1}$ is a $r$-minimal 
hypersurface of $\mathbb S^n$ if and only if the cone $C_\Sigma$ is a $r$-minimal 
hypersurface in $\mathbb R^{n+1}.$
It is a remarkable fact that a Simons' type theorem still holds for $r$-minimal cones. 
It was proved in appropriated dimensions by 
Barbosa and do Carmo in \cite{Barbosa_doCarmo} for zero scalar cones with $S_3\neq 0$ 
and by Barros and Sousa in \cite{BarrosSousa} for $r$-minimal cones assuming that 
$S_{r+2}$ is a non-null constant. 
In this setting we have the following theorem:
%%%%%%%%%%%%%%%%%%%%%%%%%%%%%%%%%%%%%%
\begin{theorem}\label{main2}
Let $C_\Sigma \subset \mathbb R ^{n+1}$, $n\geq 3$, be an unstable $r$-minimal cone for 
some  $r\in\{2, \ldots, n-1\}$. Assume that $S_j>0$, $j=1, \ldots, r-1$.
Then there is no $r$-minimal hypersurface $M^n$ properly immersed in $\mathbb R^{n+1}\setminus C_\Sigma$ such that $S_{r+1}\neq 0$.
\end{theorem}
%%%%%%%%%%%%%%%%%%%%%%%%%%%%%%%%%%%%%%

The proof follows the same steps as in the case of usual minimal hypersurfaces. In fact, 
under the hypotheses we are assuming, the linearized operator associated to the
second variation of the $r$-area functional is elliptic and the Maximum Principle  holds
(see \cite{HL, HL2}).

\subsection*{Acknowledgments}
The authors are grateful to 
Fernando Cod\'a Marques, Luciano Mari, Marcio Batista and Beno\^it Daniel 
for helpful conversations and valuable comments that have improved this work. 
The first author is partially supported by 
CNPq-Brazil (Grant 309543/2015-0), 
CAPES-Brazil (Grant 897/18) and 
FAPEAL-Brazil (Edital 14/2017).

\section{Preliminaries}
In this section we fix  notations and recall some facts we need in our proofs.

\subsection{Stability of minimal hypersurfaces}  

Let $M^n$ be an oriented Riemannian manifold and let $x: M^n\to \mathbb R^{n+1}$ be
a minimal immersion. 
If $X:M\times (-\varepsilon, \varepsilon)\to \mathbb R^{n+1}$ 
is a  smooth variation of $M$ by  immersions, whose variational vector field is $\varphi N$, 
where $N$ is a unit normal vector field along $M$ and $\varphi\in C^\infty_0(M)$, then
\[
 \frac{d \mathcal A }{dt}(t) = -n\int_M H(t)\varphi\, dM,
 \]
where $\mathcal A(t) = \int_MdM_t$ is the area functional  and
$H(t)$ is the mean curvature of the variation $X(.,t)$. 
At $t=0$ the second variation read as
\[
 \frac{d^2 \mathcal A }{dt^2}(0) =%\mathcal I(f,f)=
\int_M\left(\|\nabla \varphi\|^2 -\|A\|^2\varphi^2 \right)dM = -\int_M \varphi J\varphi\, dM.
\]
In particular, we have that $nH'(0) = J\varphi$.
In the second variation formula above, $A$ denotes the shape operator of $M$ with respect to $N$
and $J=\Delta +\|A\|^2$ is the so called \emph{Jacobi operator} or 
\emph{stability operator} of $M$.

We say that $M$ is \emph{stable} iff $\frac{d^2 \mathcal A }{dt^2}(0)\geq 0$, 
for all $\varphi\in C^\infty_0(M)$.
It is equivalent to say that the first eigenvalue of the Dirichlet problem to $-J$ 
is positive for all compact domain with regular boundary $D$ in $M$ (see \cite{FS}).

\subsection{Minimal hypersurfaces vanishing higher order mean curvature}
\label{rminimal}

Let $S_r = \sum_{i_1<\ldots <i_r} k_{i_1}\ldots k_{i_r}$, $1\leq i\leq n $, be the elementary 
symmetric functions of the principal curvatures $k_1, \ldots, k_n$ of $x: M^n\to \mathbb R^{n+1}$.
The $r$-mean curvatures of $M$ are  defined as
\[
H_0=1, \, H_r={n \choose r}^{-1}S_r,  \, r=1, \dots, n.
\]

It is well known (see \cite{Reilly}) that $r$-minimal hypersurfaces, that is $H_{r}=0$, are
critical points of $\mathcal A_r= \int_M S_{r-1}\,dM$, for all supported compact variations.
Note that the case $r=1$ corresponds to the classical case of minimal hypersurfaces.
It turns out that, in this general case, the linearized operator of $S_{r}$ is given by 
(see \cite{Rosenberg} and \cite{Reilly})
$L_r \varphi  = \textrm{div}P_r \nabla \varphi$, 
where $P_r$ are the  Newton operators, which 
are defined inductively by 
\[
P_0=Id, \quad \textrm{ and }\quad P_r=S_rId -AP_{r-1}.
\]

In this context, the second variation formula to $\mathcal A_r$ is given by
\[
\frac{d^2 \mathcal A _r}{dt^2}(0)  = -\int_M \varphi J_r\varphi\, dM,
\]
where $J_r=L_r -(r+2)S_{r+2}$. Analogously to the standard case, 
we say that a $r$-minimal hypersurface is $r$-\emph{stable} 
if  $\frac{d^2 \mathcal A_r }{dt^2}(0)\geq 0$,  for all $f\in C^\infty_0(M)$. 
It was proved by Hounie and Leite (\cite{HL,HL2}) that $L_r$ is elliptic if 
and only if $S_{r+2}\neq 0$.
So, under this condition we have the  characterization of $r$-stables
hypersurfaces in terms of the first eigenvalue of $L_r$ as in the classical case $r=1$.

We conclude this section recalling the tangency principle for $r$-minimal hypersurfaces
proved by Hounie and Leite (see \cite[Theorem  1.3 (a)]{HL2}). Namely

\begin{theorem} \label{princ}
Let $M_1$ and $M_2$ be oriented hypersurfaces in $\mathbb R^{n+1}$
satisfying 
\[
H_{r}(M_1)\leq 0 \leq  H_{r}(M_2),
\]
for a given $2\leq r< n$. Assume that $M_1$ and $M_2$ have the same normal vectors
at a tangency point $p$, with $H_j(M_2)\geq 0$, for all $j=1,\ldots, r-1$. 
If either $M_1$ or $M_2$ satisfies $H_{r+2}\neq 0$, then $M_1$ cannot remain 
above $M_2$ in a neighborhood of $p$, unless the hypersurfaces coincide locally. 

\end{theorem}

\section{Proof of the Theorems}

We start presenting below the proof of Theorem \ref{main}. 
So, assume that such $M$ does exist and let $\ell = dist(M, C_\Sigma)$ be the Euclidean distance between $M$ and $C_\Sigma$.

We first note that the distance cannot be achieved at finite points of $M$ and $C_\Sigma$.
In fact, if we assume that there exist $p\in M$ and $q\in C_\Sigma$ such that 
$\ell = dist(p,q)$ with $q\neq 0$, we can translate $M$ by $V=p-q$ until touching  
$C_\Sigma$ at  a first interior contact point and this contradicts by the Maximum Principle.

Now assume that  $q=0.$ 
In this case, we argue as 
\cite[Theorem 5.2]{MarquesNeves}.
Take $\gamma:[0,\ell]\to \mathbb R^{n+1}$ the line segment joining 
$\gamma(0)=0\in C_\Sigma$ and $\gamma(\ell)=p\in M$. 
Given $r\in (0,\ell/2)$ we have that the open ball 
$B(\gamma(r),r)$ does not intersect $C_\Sigma$, 
because $\gamma$ is the short path between $M$ and $C_\Sigma$. It implies
that the cone $C_\Sigma$ is contained in a halfspace, forcing $\Sigma $ to be totally geodesic.
A contraction.

The distance between $M$ and $C_\Sigma$ is then achieved at a sequence of points. 
Let $\{p_i\}_i\subset M$  and $\{q_i\}_i\subset C_\Sigma$ be such that
$ \ell = \lim dist(p_i,q_i). $
Note that $\{p_i\}_i$ and $\{q_i\}_i$ diverge to infinity. 
Let denote by 
\[
U_i = \frac{p_i-q_i}{\|p_i-q_i\|}\in \mathbb S^n
\]
the sequence of normalized directions. So, up a subsequence, we have $U_i\to U_0\in \mathbb S^n.$
%
%
%
%
%Writing $q_i=t_ix_i$, with $t_i\in \mathbb R_+$ and $x_i\in \Sigma$ we have that, 
%up a subsequence, $x_i\to x_0$ and so 
%\[
%\ell = \lim dist(p_i,t_ix_0).
%\]
%Let $N\in \mathbb R^{n+1}$ be a  unit parallel vector field normal to $C_\Sigma$ 
%along the generatrix $\mathbb R_+x_0$. 
Translating $M$ by $-\ell U_0$  we can assume that $dist(M, C_\Sigma)=0$. 

Now, choose $\vep\in(0,1)$ such that the $\vep$-truncated cone 
$C_\Sigma^\varepsilon = [\varepsilon,1]\times\Sigma$
is unstable (c.f. \cite[Lemma 6.1.1]{Simons}). Thus, there exist $\lambda<0$ and
a  smooth function
$\varphi: [\vep,1]\times \Sigma \to\mathbb R$ such that 
\begin{equation}\label{eq}
\begin{cases}
-J\varphi = \lambda \varphi \quad{\textrm{ in }}  C_\Sigma^\vep, \\
\varphi(\vep,.)=\varphi(1,.)=0,
\end{cases}
\end{equation}
where $J = \Delta + \|A\|^2$ is the Jacobi operator of $C_\Sigma$. 
Choosing $\lambda=\lambda_1$, the first eigenvalue of $J$, we may assume that 
$\varphi>0$ on $(\vep,1)\times\Sigma$. 
We will use $\varphi$ to construct  a perturbation of $C_\Sigma$. 
Given $\tau>0$ small and  $s\in (-\tau,\tau)$ we set
\[
D(s) = \{p+s\varphi(p)N(p):\ p\in C_\Sigma^\vep \},
\]
where $N$ is the unit normal vector on $C_\Sigma$ pointing inward $C_\Sigma^+$, 
the  connected component of $\mathbb R^{n+1}\setminus C_\Sigma$ that contains 
$M$.
Since $M$ is proper, there exists $\tau_0\leq \tau$ such that
$D(s)\cap M = \emptyset$, for all $s\in(\tau_0,\tau_0)$. 
%Moreover, the mean curvature vector $\vec H(s)$ of $D(s)$ is   

We fix $s\in(-\tau_0, \tau_0)$, $s\neq 0$, such that $D(s)$ is contained in $C_\Sigma^+$
and let us denote by $H(s)$ the mean curvature function of 
$D(s)$ with respect to the unit normal $N(s)$ pointing inward $C_\Sigma^+$.
Since $\varphi$ is a non trivial solution of problem (\ref{eq}) we have that 
\[
H'(0)=\frac 1 nJ\varphi = -\lambda \varphi>0
\]
and thus $H(s)>0 $ for the desired $s$.

Now, we  consider the perturbed cone defined as 
\[
\widetilde C_\Sigma=
\begin{cases}
D(s), \textrm{ if } t\in(\vep, 1)\\
C_\Sigma, \textrm{ othewise}.
\end{cases}
\]
 
Since $M$ is asymptotic to $C_\Sigma$ if we apply a family of homotheties $R_a$ to $M$
we can find $a\in (0,1)$ such that $M_a=R_a(M)$ is tangential to $\widetilde C_\Sigma$ at
some point of $D(s)$. Since homothety preserves the minimality 
we have a contradiction with the Maximum Principle. 

\medskip
The proof of Theorem \ref{main2} follows
the same steps as above. In fact, homotheties preserve $r$-minimality and
we can choose the perturbations of the cone too small
such that the conditions on the positivity of $S_j$ are preserved and so the
 tangency principle (Theorem \ref{princ}) applies.

%\qed

\bibliographystyle{amsplain}
\bibliography{bibliography}

\end{document}